\documentclass[12pt]{article}
\usepackage {amsmath}
\input amssym.def
\input amssym.tex

\def \C{{\Bbb C}}
\def \Q{{\Bbb Q}}

\def \proof{{\noindent{\it Proof.\ \ }}}
\newtheorem{Th}{THEOREM}[section]

\newtheorem{prop}[Th]{PROPOSITION}

\newtheorem{lemma}[Th]{LEMMA}

\title{Equivariant Singular Riemann-Roch Theorem}
\author{Bin Zhang} 
\begin{document}
\maketitle
\section {Introduction}

The Rieman-Roch theorem \cite {hf1}, \cite {bfm1} has played
 such a key role in the study of algebraic geometry, algebraic 
 topology or even number theory, that many attempts have been made
 to generalize it  to the equivariant case (for example \cite {kt1}, \cite {bv1}). 
 Here we consider the case of the equivariant singular Riemann-Roch theorem, 
 following the framework of P. Baum, W. Fulton and  R. MacPherson.

\section {Totaro's approximation of $EG$}

The purpose of this section is to describe a construction, due to Totaro \cite {eg1}, of an algebro-geometric substitute for the classifying space of topological group.

For topological group $G$, the spaces $EG$ and $BG$ are infinite dimensional, it is hard to control in general. What Totaro produces for a reductive algebraic group is a directed system of $G$-bundles $E_n \to B_n$ with the following property:
for any algebraic principal $G$-bundle $E\to X$, there is a map $X'\to X$ with the fiber isomorphic to affine space ${\Bbb A}^m$, such that the pullback bundle $E'\to X'$ is pulled back from one of the bundles in the directed system by a map $X'\to B_n$.

For a complex linear reductive algebraic group $G$, Totaro constructs the directed system as follows:

Every object in the directed system is a pair $(V, V')$, such that $V$ is a representation of $G$, $V'$ is a non-empty open set of $V$, $G$ acts freely on $V'$, $V'\to V'/G$ is a principal $G$-bundle. $(V, V')<(W, W')$ if $V$ is a subspace of $W$, $V'\subset W'$, ${\rm codim}_V(V-V')<{\rm codim}_W(W-W')$. The morphisms are just inclusions.

For some groups, there is a convenient choice of the directed system. For example, if $G=T$ is a split torus of rank n, then $\{(V^n, (V-\{0\})^n)\ | \ {\rm dim} V=l\}$ is the directed system.

\section {Equivariant cohomology}

In general \cite {ma1}, we define 

$$H_G^*(X)=H^*(X\times ^GEG),$$

Using the directed system above, 

\begin {lemma} For a complex Lie group $G$ and for $X$ a complex algebraic variety on which $G$ operates algebraically, we have
$$H^*_G(X)=\lim _{\underset {(V, V')} \gets} H^*(X\times ^G V')$$
\end {lemma}

\proof See \cite {bm1} $\blacksquare$

For every $G$-vector bundle $E$ over $X$, $E\times ^G V'\to X\times ^G V'$ and $E\times ^G EG\to X\times ^G EG$ are all $G$-vector bundles, so we can define the Chern character $ch^G(E) \in H^*_G(X)$ and $ch^V(E) \in H^*(X\times ^G V')$by 
$$ch^G(E)=ch (E\times ^G EG)$$
and
$$ch^V(E)=ch (E\times ^G V')$$
We know that under the map $H^*_G(X) \to H^*(X\times ^G V')$, $ch^G(E)$ goes to $ch^V(E)$.
 
When $X$ is smooth, then for any $(V, V')$ in the directed system, $X\times ^G V'$ is smooth. So we can define the cohomological equivariant Todd class $Td^G(X)$ to be the element in $H_G^*(X)$ which maps to $Td(X\times ^GV')$ for every $(V, V')$ in the directed system.
  
\section {Equivariant homology}

Now for a complex linear reductive algebraic group $G$, by using the constructed directed system, we can define the equivariant homology \cite {eg1} for a $G$-variety $X$, suppose dim$_{\C}X$=n.

$$H_i^G(X)=\lim _{\underset {(V, V')} \gets} H^{BM}_{i+2l-2g}(X\times ^GV')$$

where $l={\rm dim}_{\C}V$, $g={\rm dim}_{\C}G$.

{\bf Remark} It is possible that $H_i^G(X) \not =0$ for negative $i$.

\begin {prop} $H^{G}_*(X)$ is independent of $(V, V')$ if codim$_V (V-V')$ big enough. 
\end {prop}

$H^G_*(X)$ is a module over $H_G^*(X)$. In other words there is a cap product:
$\cap:  H_G^*(X)\otimes H^G_*(X)\to H^G_*(X)$ defined in the following way:

If $\alpha \in H_G^*(X)$, 
then for any $(V,V')$, the image of $\alpha$ under the map $H_G^*(X)\to H^*(X\times ^G V')$ defines a map by cap product: $H^{BM}_*(X\times ^G V') \to H^{BM}_*(X\times ^G V') \to H^{G}_*(X)$. Go over to the direct limit, these 
maps induce a map from $H^{G}_*(X)$ to $H^{G}_*(X)$. This gives us the module structure.

For a non-singular variety $X$, we can define the orientation class $[X]_G \in H^{G}_{2n}(X)$ as follows:

For any object $(V, V')$ in the directed system, $X\times ^G V'$ is non-singular, so we have a class $[X\times ^G V'] \in H^{BM}_{2n+2l-2g}(X\times ^GV')$, $\{ [X\times ^G V']\}$ gives the orientation class $[X]_G$.

\section {Equivariant Riemann-Roch Theorem}

Now we can define the equivariant Todd class $\tau ^G$ as follows,

For any equivariant coherent $G$-sheaf ${\cal F}$ over $X$ ($X$ is a complex algebraic variety), if $p: \ X\times V' \to X$ is the projection, then the pullback sheaf $p^*{\cal F}$ is a sheaf over $X\times V'$, it descents to a coherent sheaf ${\cal F}_V$ over $X\times ^G V'$. $\tau ({\cal F}_V)$ is in $H^{BM}_*(X\times ^G V')$, so we get a map $K^G_0(X)\to H^{BM}_*(X\times ^G V')$. Obviously, it is compatible with the transition maps, so we get a map  $K^G_0(X)\to H^{G}_*(X)$, that is $\tau ^G$.

\begin {Th} [Equivariant Riemann-Roch Theorem]
 $\tau ^G$ is the unique natural transformation (with respect to $G$-equivariant proper maps) from $K_0^G(X)$ to $H^G_*(X)$, such that for any $G$-variety $X$, we have the following commutative diagram:

$$\begin {array}{cclcc} 
K_G^0(X)&\otimes& K_0^G(X) &{\overset {\otimes} \to}& K_0^G(X)\\
&\downarrow&ch^G\otimes \tau ^G&&\downarrow\tau ^G\\
H_G^*(X)&\otimes& H_*^G(X) &{\overset {\cap} \to} &H_*^G(X)
\end {array}
$$
and if $f: X\to Y$ is a proper $G$-equivariant algebraic map, then 
$$\begin {array}{rrr} 
K_0^G(X) &{\overset {f_*} \to}& K_0^G(Y)\\
\downarrow \tau ^G&&\downarrow\tau ^G\\
H_*^G(X) &{\overset {f_*} \to} &H_*^G(X)
\end {array}
$$
Furthermore, if $X$ is non-singular, and ${\cal O}_X$ is the structure sheaf, then 
$$
\tau ^G({\cal O}_X)=Td^G(X) \cap [X]_G
$$
\end {Th}
\proof In fact, the proof is already encoded in the definition. $\blacksquare$

From the discussion, we have the following commutative diagram:

$$\begin {array}{ccc} 
K_0^G(X) &{\overset {\tau ^G} \to}& H_*^G(X)\\
\downarrow&&\downarrow\\
K_0(X) &{\overset {\tau} \to} &H_*(X)
\end {array}
$$
where both downarrows are induced by forgetting functors.

\section {Equivariant Riemann-Roch theorem with value in equivariant Chow group}

For linear algebraic group $G$, by using the same directed system, we can define \cite {eg2}

$$A_i^G(X)={\underset {\to} \lim} A_{i+2l-2g}(X\times ^GV')$$
where $l={\rm dim}_{\C}V$, $g={\rm dim}_{\C}G$.

For the equivariant operational Chow group $A^*_G(X)$, we define it as
$$A^i_G(X)={\underset {\gets} \lim} A^{i}(X\times ^GV')$$
This definition concises with the definition in \cite {eg2}.

As in the previous section, we can define the Chern character, the Todd class and get a theorem similar to Theorem 5.1, except we use $A^i_G(X)\otimes {\Q}$ and$A_*^G(X)\otimes {\Q}$ to replace $H^*_G(X)$ and $H_*^G(X)$.

\end{document}